\documentclass[12pt]{article}
\usepackage{graphicx}
\usepackage{amsmath,amsthm,amssymb,enumerate}
\usepackage{euscript,mathrsfs}
\usepackage{color}
\usepackage{dsfont}
\usepackage[left=2cm,right=2cm,top=3.5cm,bottom=3.5cm]{geometry}
\usepackage{color}
\usepackage[framemethod=tikz]{mdframed}
\allowdisplaybreaks

\usepackage{soul}

\catcode`\@=11 \@addtoreset{equation}{section}

\catcode`\@=12

\newtheorem{Theorem}{Theorem}[section]
\newtheorem{Proposition}[Theorem]{Proposition}
\newtheorem{Lemma}[Theorem]{Lemma}
\newtheorem{Corollary}[Theorem]{Corollary}

\theoremstyle{definition}
\newtheorem{Definition}[Theorem]{Definition}

\newtheorem{Remark}[Theorem]{Remark}

\newcommand{\bTheorem}[1]{
	\begin{Theorem} \label{T#1} }
	\newcommand{\eT}{\end{Theorem}}

\newcommand{\bProposition}[1]{
	\begin{Proposition} \label{P#1}}
	\newcommand{\eP}{\end{Proposition}}

\newcommand{\bLemma}[1]{
	\begin{Lemma} \label{L#1} }
	\newcommand{\eL}{\end{Lemma}}

\newcommand{\bCorollary}[1]{
	\begin{Corollary} \label{C#1} }
	\newcommand{\eC}{\end{Corollary}}

\newcommand{\bRemark}[1]{
	\begin{Remark} \label{R#1} }
	\newcommand{\eR}{\end{Remark}}

\newcommand{\bDefinition}[1]{
	\begin{Definition} \label{D#1} }
	\newcommand{\eD}{\end{Definition}}

\newcommand{\Del}{\Delta_x}

\newcommand{\Ds}{\mathbb{D}_x}

\newcommand{\vuB}{\vc{u}_B}

\newcommand{\data}{{\rm data}}

\newcommand{\bfphi}{\boldsymbol{\varphi}}

\newcommand{\bFormula}[1]{
	\begin{equation} \label{#1}}
	\newcommand{\eF}{\end{equation}}

\newcommand{\Ov}[1]{\overline{#1}}

\newcommand{\aleq}{\stackrel{<}{\sim}}

\newcommand{\vr}{\varrho}
\newcommand{\vre}{\vr_\ep}

\newcommand{\vte}{\vt_\ep}
\newcommand{\vue}{\vu_\ep}

\newcommand{\tvt}{\tilde \vt}

\newcommand{\vt}{\vartheta}
\newcommand{\vu}{\vc{u}}

\newcommand{\vc}[1]{{\bf #1}}

\newcommand{\Div}{{\rm div}_x}
\newcommand{\Grad}{\nabla_x}

\newcommand{\dx}{\,{\rm d} {x}}

\newcommand{\dt}{\,{\rm d} t }

\newcommand{\dxdt}{\dx  \dt}
\newcommand{\intO}[1]{\int_{\Omega} #1 \ \dx}

\newcommand{\intST}[1]{\int_{S_T} #1 \ \dt}

\newcommand{\D}{{\rm d}}

\newcommand{\ep}{\varepsilon}

\newcommand{\vtB}{\vt_B}

\newcommand{\br}{ \nonumber \\ }

\def\softd{{\leavevmode\setbox1=\hbox{d}%
		\hbox to 1.05\wd1{d\kern-0.4ex{\char039}\hss}}}
\definecolor{Cgrey}{rgb}{0.85,0.85,0.85}
\definecolor{Cblue}{rgb}{0.50,0.85,0.85}
\definecolor{Cred}{rgb}{1,0,0}
\definecolor{fancy}{rgb}{0.10,0.85,0.10}

\newcommand\Cbox[2]{%
	\newbox\contentbox%
	\newbox\bkgdbox%
	\setbox\contentbox\hbox to \hsize{%
		\vtop{
			\kern\columnsep
			\hbox to \hsize{%
				\kern\columnsep%
				\advance\hsize by -2\columnsep%
				\setlength{\textwidth}{\hsize}%
				\vbox{
					\parskip=\baselineskip
					\parindent=0bp
					#2
				}%
				\kern\columnsep%
			}%
			\kern\columnsep%
		}%
	}%
	\setbox\bkgdbox\vbox{
		\color{#1}
		\hrule width  \wd\contentbox %
		height \ht\contentbox %
		depth  \dp\contentbox
		\color{black}
	}%
	\wd\bkgdbox=0bp%
	\vbox{\hbox to \hsize{\box\bkgdbox\box\contentbox}}%
	\vskip\baselineskip%
}

\mdfdefinestyle{MyFrame}{%
	linecolor=black,
	outerlinewidth=1pt,
	roundcorner=5pt,
	innertopmargin=\baselineskip,
	innerbottommargin=\baselineskip,
	innerrightmargin=10pt,
	innerleftmargin=10pt,
	backgroundcolor=white!20!white}

\date{}

\begin{document}


\title{Taylor--Couette flow with temperature fluctuations: Time periodic solutions}

\author{Eduard Feireisl
	\thanks{The work of E. Feireisl was partially supported by the
		Czech Sciences Foundation (GA\v CR), Grant Agreement
		21--02411S. The Institute of Mathematics of the Academy of Sciences of
		the Czech Republic is supported by RVO:67985840.} \and Young--Sam Kwon
	\thanks{The work of Y.--S. Kwon was partially supported by
		the National Research Foundation of Korea (NRF2020R1F1A1A01049805)}}

\date{\today}

\maketitle

\bigskip

\centerline{Institute of Mathematics of the Academy of Sciences of the Czech Republic}

\centerline{\v Zitn\' a 25, CZ-115 67 Praha 1, Czech Republic}

\medskip

\centerline{Department of Mathematics, Dong-A University}

\centerline{Busan 49315, Republic of Korea}

\bigskip

\begin{abstract}
	
	We consider the motion of a viscous compressible and heat conducting fluid confined in the gap between
	two rotating cylinders (Taylor--Couette flow). The temperature of the cylinders is fixed but not necessarily
	constant. We show that the problem admits a time--periodic solution as soon as the ratio of the angular velocities of the two cylinders is a rational number.
	
\end{abstract}




\tableofcontents

\section{Introduction}
\label{i}

The Taylor--Couette flow is one of the iconic examples of turbulent fluid motion driven
by the adherence of a viscous fluid to the kinematic boundary. Specifically, consider two
not necessarily concentric cylinders
\begin{align}
\mathcal{C}_1 = B_{\underline{r}} (\vc{a}) \times [0,1]&,\ \mathcal{C}_2 = B_{\Ov{r}}(\vc{b}) \times [0,1],\ \br
B_{\underline{r}} (\vc{a}) = \left\{ x_h \equiv (x_1,x_2) \ \Big| \ |x_h - \vc{a}| \leq \underline{r} \right\}&,
\ B_{\Ov{r}} (\vc{b}) = \left\{ x_h \equiv (x_1,x_2) \ \Big| \ |x_h - \vc{b}| \leq \Ov{r} \right\}, \br
0 < \underline{r}&, \ |\vc{a} - \vc{b}| + \underline{r} < \Ov{r}. \nonumber
\end{align}
For the sake of simplicity, we suppose the motion is periodic in the vertical direction $x_3$ therefore the
fluid domain $\Omega \subset R^3$ is given as
\begin{equation} \label{i1}
	\Omega = \Big( {\rm int} [ \mathcal{C}_2 ] \setminus \mathcal{C}_1 \Big) \times \mathbb{T}^1,\
	\mathbb{T}^1 \equiv [0,1] \Big|_{\{ 0,1 \} }.
	\end{equation}
The cylinders $\mathcal{C}_1$, $\mathcal{C}_2$ rotate with constant angular velocities $\omega_1$, $\omega_2$ respectively. Finally, we suppose the surface temperatures  $\theta_1$ of $\partial \mathcal{C}_1$,
$\theta_2$ of $\partial \mathcal{C}_2$ are independent of time but not necessarily constant.

We consider a general viscous and heat conducting fluid. Accordingly, the fluid velocity $\vu = \vu(t,x)$
satisfies the no--slip boundary conditions
\begin{equation} \label{i2}
	\vu|_{\partial \Omega} = \vuB, \ \mbox{where} \
	\vuB (x) = \left\{ \begin{array}{l} \omega_1 (x_h - \vc{a})^\perp \ \mbox{on} \ \partial \mathcal{C}_1,\\
	\omega_2 (x_h - \vc{b})^\perp \ \mbox{on} \ \partial \mathcal{C}_2. \end{array} \right.	
	\end{equation}

Similarly, we impose the Dirichlet boundary conditions for the fluid temperature $\vt = \vt(t,x)$. Although the temperature of the cylinders is independent of time, the fact that they rotate enforces the \emph{time dependent}
boundary conditions:
\begin{equation} \label{i3}
\vt|_{\partial \Omega} = \vtB,\ \vtB(t,x) = \left\{ \begin{array}{l} \theta_1 \left(\vc{a} + \underline{r} [\cos (\omega_ 1 t), \sin(\omega_1 t)], x_3 \right)  \ \mbox{on} \ \partial \mathcal{C}_1,
\\ \theta_2  \left(\vc{b} + \Ov{r} [\cos (\omega_ 2 t), \sin(\omega_2 t)], x_3  \right) \
\mbox{on} \ \partial \mathcal{C}_2.
 \end{array} \right.	
	\end{equation}
If the ratio
\begin{equation} \label{i3a}
	\frac{\omega_1}{\omega_2} = k - \ \mbox{a rational number}
	\end{equation}
the boundary temperature $\vt_B$ is periodic in time,
\begin{equation} \label{i3b}
	\vtB(t + T, x) = \vtB(t,x) \ \mbox{for a certain}\ T \geq 0.
\end{equation}
As pointed out in Yang et al \cite{Yangetal}: ``\emph{Many natural and industrial turbulent flows are subjected to time-dependent boundary conditions.}''

\subsection{Field equations}

The time evolution of the fluid density $\vr = \vr(t,x)$, the (absolute) temperature $\vt = \vt(t,x)$ and the velocity $\vu = \vu(t,x)$ is governed by the \emph{Navier--Stokes--Fourier (NSF) system}  of partial differential equations:

\begin{align}
	\partial_t \vr + \Div (\vr \vu) &= 0, \label{i4} \\
	\partial_t (\vr \vu) + \Div (\vr \vu \otimes \vu)  + \Grad p (\vr, \vt) &= \Div \mathbb{S}(\vt, \Ds \vu) + \vr \Grad G, \label{i5}\\
	\partial_t (\vr e(\vr, \vt)) + \Div (\vr e(\vr, \vt) \vu) + \Grad \vc{q}(\vt, \Grad \vt) &= \mathbb{S}
	(\vt, \Ds \vu) : \Ds \vu - p(\vr, \vt) \Div \vu,
	\label{i6}
	\end{align}

\noindent where $\mathbb{S}$ is the viscous stress given by Newton's rheological law
\begin{equation} \label{i7}
	\Ds \vu = \frac{1}{2} \left( \Grad \vu + \Grad^t \vu \right),\ \mathbb{S} (\vt, \Ds \vu) = \mu(\vt) \left( \Grad \vu + \Grad^t \vu - \frac{2}{3} \Div \vu \mathbb{I} \right) +
	\eta(\vt) \Div \vu \mathbb{I},
\end{equation}
and $\vc{q}$ is the heat flux given by
Fourier's law
\begin{equation} \label{i8}
	\vc{q}(\vt, \Grad \vt) = - \kappa(\vt) \Grad \vt.
\end{equation}
The pressure $p=p(\vr, \vt)$ and the internal energy are interrelated by Gibbs' equation
\begin{equation}\label{i8b}
\vt Ds = De + p D \frac{1}{\vr},
\end{equation}
where $s$ is a new thermodynamic function called entropy.

\subsection{Time periodic solutions}

Motivated by the example of the Taylor--Couette flow, our goal is to establish existence of time--periodic solutions
to the NSF system driven by time periodic boundary conditions. Specifically, we consider a bounded domain $\Omega \subset R^3$,
\begin{align}
\Omega \ \mbox{of class}\ C^\infty,\ \partial \Omega &= \cup_{i=1}^n \Gamma_i, \br
\vuB &= \vuB(x), \vuB \cdot \vc{n} |_{\partial \Omega} = 0, \br
\vtB (t + T, \cdot)|_{\partial \Omega} &= \vtB(t,\cdot) |_{\partial \Omega} \ \mbox{for some}\ T > 0.
\label{dom}
\end{align}
The assumption $\vuB$ independent of time can be replaced by time--periodicity, the assumption that $\vuB$ is tangential to the boundary is however essential as it entails the total mass conservation
\begin{equation} \label{tmc}
M \equiv	\intO{ \vr(t, \cdot) } \ \mbox{for any}\ t.
\end{equation}

The present work can be seen as a continuation of the papers \cite{FeGwSG22}
and \cite{FeMuNoPo}  devoted to the
time periodic solutions to the NSF system. The main difference between \cite{FeGwSG22} and \cite{FeMuNoPo}
is the concept of weak solution. The approach of \cite{FeMuNoPo} is based on the mathematical theory developed
in the monograph \cite{FeNo6A} for essentially energetically closed fluid systems, while the more recent
result \cite{FeGwSG22} requires the new concept of weak solution for general open systems introduced in
\cite{ChauFei} and further elaborated in the monograph \cite{FeiNovOpen}. Similarly, the present paper
needs the new framework \cite{FeiNovOpen}.

At first glance, the results presented below could be seen as a generalization of \cite{FeGwSG22} to the case of inhomogeneous boundary velocity. There is, however, a substantial difference due to the choice of the equation of state. In order to handle the problem driven by the motion of the boundary, the pressure equation of state must be augmented by the so--called hard--sphere pressure component already used in \cite{FeiKwo}, see also \cite[Chapter 9, Section 9.1.4]{FeiNovOpen}. Such a hypothesis is not needed in \cite{FeGwSG22} therefore the principal part of the analysis based on {\it a priori} bounds in \cite{FeGwSG22} is quite complementary to the present case.

To conclude the introductory part, let us recall that there several results concerning the time periodic
solutions to the compressible and/or heat conducting fluid systems driven by \emph{smooth and small} data, see
e.g. B\v rezina and Kagei \cite{BreKag2}, \cite{BreKag1}, Jin and Yang \cite{JinYang} , Kagei and Oomachi \cite{KagOom},
Kagei and Tsuda \cite{KagTsu}, Tsuda \cite{Tsuda}, Valli and Zajaczkowski \cite{Valli1}, \cite{VAZA} to name only a few. Last but not least, it is worth mentioning there are alternative approaches to the concept of weak solutions 
proposed Bresch and Desjardins \cite{BRDE}, \cite{BRDE1} or Bresch and Jabin \cite{BreJab}. Neither of them, however, seem to apply to problems with inhomogeneous boundary conditions.

The paper is organized as follows. In Section \ref{m}, we list the main hypotheses and state our main result.
In Section \ref{A}, we introduce a family of approximate problems solvable by the methods developed in
\cite{FeMuNoPo}. In Section \ref{e}, we perform the limit in the sequence of approximate solutions and prove the main existence result. Section \ref{CR} contains concluding remarks.

\section{Hypotheses and main result}
\label{m}

In this section, we collect the necessary hypotheses imposed on the constitutive equations and state our main result.

\subsection{Hypotheses}

Following \cite{FeiKwo}, \cite[Chapter 4, Section 4.3]{FeiNovOpen} we impose the following hypotheses concerning
the equation of state:
\begin{equation} \label{m1}
	p(\vr ,\vt) = p_m (\vr, \vt) + p_r(\vt) + p_{HS}(\vr),\
	e(\vr, \vt) = e_M (\vr, \vt) + e_r(\vr, \vt) + e_{HS} (\vr),
	\end{equation}
where
\begin{align}
	p_{m}(\vr, \vt)  &= \vt^{\frac{5}{2}} P \left( \frac{\vr}{\vt^{\frac{3}{2}}  } \right),\
	p_r (\vt) = \frac{a}{3} \vt^4,\ a > 0,\  p_{HS}(\vr) = \frac{b \vr}{(\Ov{\vr} - \vr)^s },\ b > 0,
	\label{m2} \\
	e_m (\vr, \vt) &= \frac{3}{2} \frac{\vt^{\frac{5}{2}} }{\vr} P \left( \frac{\vr}{\vt^{\frac{3}{2}}  } \right),\ e_r (\vr, \vt) =  \frac{a}{\vr} \vt^4,\  e_{HS}(\vr) = \int_{\Ov{\vr}/2}^\vr \frac{b}{z (\Ov{\vr} - z)^s } \D z,
	\ s > 3,
	\label{m3}
\end{align}
where $P \in C^1[0,\infty)$ satisfies
\begin{equation} \label{m4}
	P(0) = 0,\ P'(Z) > 0 \ \mbox{for}\ Z \geq 0,\ 0 < \frac{ \frac{5}{3} P(Z) - P'(Z) Z }{Z} \leq c \ \mbox{for}\ Z > 0.
\end{equation} 	
In particular, the function $Z \mapsto P(Z)/ Z^{\frac{5}{3}}$ is decreasing, and we suppose
\begin{equation} \label{m5}
	\lim_{Z \to \infty} \frac{ P(Z) }{Z^{\frac{5}{3}}} = p_\infty \geq 0.
\end{equation}
The associated  entropy $s$ reads
\begin{equation} \label{m6}
	s(\vr, \vt) = s_m (\vr, \vt) + s_r(\vr, \vt),\ s_m(\vr, \vt) = \mathcal{S} \left( \frac{\vr}{\vt^{\frac{3}{2}} } \right),\ s_r(\vr, \vt) = \frac{4a}{3} \frac{\vt^3}{\vr},
\end{equation}
where
\begin{equation} \label{m7}
	\mathcal{S}'(Z) = -\frac{3}{2} \frac{ \frac{5}{3} P(Z) - P'(Z) Z }{Z^2}.
\end{equation}
Finally, in accordance with the Third law of thermodynamics, we suppose
\begin{equation} \label{m8}
	\lim_{\vt \to 0} s(\vr, \vt) = 0 \ \mbox{for any fixed}\ \vr > 0, \ \mbox{meaning}\
	\mathcal{S}(Z) \to 0 \ \mbox{as}\ Z \to \infty.
\end{equation}
The reader may consult \cite[Chapter 4]{FeiNovOpen} for the physical background of the above hypotheses.

In addition,
the transport coefficients $\mu$, $\eta$, and $\kappa$ are continuously differentiable functions of the temperature $\vt$ satisfying
\begin{align}
	0 < \underline{\mu} \left(1 + \vt \right) &\leq \mu(\vt) \leq \Ov{\mu} \left( 1 + \vt \right),\
	|\mu'(\vt)| \leq c \ \mbox{for all}\ \vt \geq 0,\ \frac{1}{2} \leq \Lambda \leq 1, \br
	0 &\leq  \eta(\vt) \leq \Ov{\eta} \left( 1 + \vt \right), \br
	0 < \underline{\kappa} \left(1 + \vt^\beta \right) &\leq  \kappa(\vt) \leq \Ov{\kappa} \left( 1 + \vt^\beta \right),\ \beta > 6.
	\label{m9}
\end{align}

The pressure as well as the internal energy are augmented by the so--called hard--sphere component
$p_{HS}$, $e_{HS}$ respectively. They are both singular at $\Ov{\vr}$ and force the density
$\vr$ to be bounded above by $\Ov{\vr}$. This facilitates considerably the analysis leading to new {\it a priori}
estimates. The regularizing effect of the hard--sphere pressure has been exploited in a number of recent studies:
\cite{AbbFei21per}, \cite{CiFeJaPe1},
\cite{FeiZha}, the monograph \cite{FeiNovOpen} and the references cited therein.

\subsection{Weak solution}

It is convenient to identify the time periodic functions with distributions defined on the flat torus
\[
S_T = [0,T]|_{\{ 0, T \} }.
\]
Suppose that $\vuB$, $\vtB$ have been extended inside $\Omega$.

\begin{Definition}[{\bf Weak solution}] \label{Dws}
	
	We say that a trio $(\vr, \vt, \vu)$ is \emph{weak time periodic solution} to the NSF system \eqref{i4} - \eqref{i8b}, with the boundary conditions \eqref{dom} if the following holds:
	\begin{itemize}
	\item {\bf Integrability.}
	\begin{align}
		0 &\leq \vr < \Ov{\vr} \ \mbox{a.a. in}\ S_T \times \Omega,\
		\vr \in C_{\rm weak}(S_T; L^q(\Omega)) \ \mbox{for any}\ 1 \leq q < \infty; \br
		0 &< \vt \ \mbox{a.a. in}\ S_T \times \Omega,\ \vt \in L^\infty(S_T; L^4(\Omega)) \cap
		L^2(S_T; W^{1,2}(\Omega)),\ \vt|_{\partial \Omega} = \vtB;\br
		\vu &\in L^2(S_T; W^{1,2}(\Omega; R^3),\ \vu|_{\partial \Omega} = \vuB; \br
		\vr \vu &\in C_{\rm weak}(S_T; L^2(\Omega; R^3)).
		\nonumber
		\end{align}
	
	\item {\bf Equation of continuity.}
	\[
		\intST{ \intO{ \left[ \vr \partial_t \varphi + \vr \vu \cdot \Grad \varphi \right] } } = 0
	\]
	for any $\varphi \in C^1(S_T \times \Ov{\Omega})$.
	
	\item {\bf Momentum balance.}
	\begin{align}
	&\intST{ \intO{ \Big[ \vr \vu \cdot \partial_t \bfphi + \vr \vu \otimes \vu : \Grad \bfphi					 +
			p (\vr, \vt) \Div \bfphi \Big] } } \br &\quad = \intST{ \intO{ \Big[ \mathbb{S}(\vt, \Ds \vu) : \Ds \bfphi - \vr \Grad G \cdot \bfphi \Big] } }
		\nonumber
		\end{align}
	for any $\bfphi \in C^1_c(S_T \times \Omega; R^3)$.
	
	\item {\bf Entropy inequality.}	
	
	\begin{align}
		- &\intST{ \intO{ \left[ \vr s(\vr, \vt) \partial_t \varphi + \vr s(\vr, \vt) \vu \cdot \Grad \varphi + \frac{\vc{q} (\vt, \Grad \vt)}{\vt} \cdot
				\Grad \varphi \right] } } \br &\geq \intST{ \intO{ \frac{\varphi}{\vt} \left[ \mathbb{S}(\vt, \Ds \vu) : \Ds \vu -
				\frac{\vc{q}(\vt, \Grad \vt) \cdot \Grad \vt }{\vt} \right] } }
		\nonumber
	\end{align}
for any $\varphi \in C^1_c(S_T \times \Omega)$, $\varphi \geq 0$.

\item {\bf Ballistic energy balance.}	

\begin{align}
	- &\intST{ \partial_t \psi	\intO{ \left[ \frac{1}{2} \vr |\vu -\vuB|^2 + \vr e - \tvt \vr s \right] } }
	+ \intST{ \psi
		\intO{ \frac{\tvt}{\vt}	 \left[ \mathbb{S}: \Ds \vu - \frac{\vc{q} \cdot \Grad \vt }{\vt} \right] } }  \br
&\leq
	\intST{ \psi \intO{ \left[ \vre (\vu - \vuB) \cdot \Grad G - \vr s \vu \cdot \Grad \tvt - \frac{\vc{q}}{\vt} \cdot \Grad \tvt
			- \partial_t \tvt \vr s
			\right] } } \br
	&- \intST{ \psi \intO{ \left[ \vr \vu \otimes \vu + p(\vr, \vt) \mathbb{I} - \mathbb{S}(\vt, \Ds
			\vu) \right] : \Ds \vuB } }\br &+ \frac{1}{2} \intST{ \psi \intO{\vr \vu \cdot \Grad |\vuB|^2 } }	
	\nonumber
\end{align}	
for any $\psi \in C^1(S_T)$, $\psi \geq 0$, and any
\begin{equation}\label{class}
\tvt \in C^1(S_T \times \Ov{\Omega}), \ \tvt > 0,\ \tvt|_{\partial \Omega} = \vtB.	
	\end{equation}	
		\end{itemize}
	
	\end{Definition}

Note carefully that the above definition is the same as in \cite[Chapter 12, Definition 7]{FeiNovOpen}. As
the density $\vr$ is bounded, we can use the regularization technique of DiPerna and Lions \cite{DL} to
obtain a renormalized version of the equation of continuity
\[
	\intST{ \intO{ \left[ b(\vr) \partial_t \varphi + b(\vr) \vu \cdot \Grad \varphi + \Big(
		b(\vr) - b'(\vr) \vr \Big) \Div \vu \varphi \right] } } =0
\]
for any $\varphi \in C^1(S_T \times \Ov{\Omega} )$, and any $b \in C^1(R)$.

\subsection{Main result}

We are ready to state our main result concerning the existence of time periodic solutions.

\begin{Theorem}[{\bf Time periodic solution}] \label{MT1}
	
	Suppose that $\Omega \subset R^3$ is a bounded domain of class \eqref{dom}, where the boundary data satisfy
	\begin{equation} \label{hyp1}
		\vuB \in C^2(\partial \Omega; R^3)), \ \vtB \in C^2(S_T \times \partial \Omega),\
		\inf_{S_T \times \partial \Omega} \vtB > 0.
		\end{equation}
	Let $G \in W^{1,\infty}(\Omega)$, and let the thermodynamic functions $p$, $e$, $s$, and the transport coefficients $\mu$, $\eta$, $\kappa$ satisfy the hypotheses \eqref{m1}--\eqref{m9}, with
	\begin{equation} \label{hyp2}
		s > 3, \ \beta > 6.
		\end{equation}
	Let
	\[
	0 < M < \Ov{\vr} |\Omega|
	\]
	be given.
	
	Then the NSF system \eqref{i4}--\eqref{i8}, with the boundary conditions \ref{dom} admits
	a time periodic solution $(\vr, \vt, \vu)$ in the sense specified in Definition \ref{Dws} such that
	\[
	\intO{ \vr(t, \cdot) } = M \ \mbox{for any}\ t \in S_T.
	\]

	\end{Theorem}

The rest of the paper is devoted to the proof of Theorem \ref{MT1}.

\section{Approximate problem}
\label{A}

Following the strategy of \cite{FeGwSG22}
we consider an approximate problem replacing the boundary conditions \eqref{dom} by
\begin{equation} \label{A1}
	\vu \cdot \vc{n} = 0,\
	[\mathbb{S} \cdot \vc{n}]_{\rm tan} = - \frac{1}{\ep} (\vu - \vuB) \ \mbox{on}\ \partial \Omega,
	\end{equation}
\begin{equation} \label{A2}
	\vc{q} \cdot \vc{n} = \frac{1}{\ep} |\vt - \vtB|^{k}(\vt - \vtB),\ k = \beta - 1, \ \mbox{on}\ \partial \Omega,
\end{equation}
where $\ep > 0$ is a small parameter. Note that \eqref{A1} can be interpreted as Navier's boundary
condition with friction while \eqref{A2} is a nonlinear Robin type boundary condition introduced in
\cite{BaFeLMMiYu}, \cite{FeGwSG22}.

Moreover, we replace the hard--sphere component of the pressure by a suitable cut--off,
Specifically, we consider
\begin{align}
p_N (\vr, \vt) &= p_m (\vr, \vt) + p_r (\vr , \vt) + \frac{1}{N} \vr^2 +
([ \vr - \Ov{\vr} - 1]^+)^\Gamma + p^N_{HS}(\vr),\br
p^N_{HS}(\vr) &= \left\{ \begin{array}{l} p_{HS} (\vr) \ \mbox{if}\ 0 \leq \vr \leq \Ov{\vr} - \frac{1}{N}, \\  \\
a_1 \vr + a_2 \ \mbox{if}\ \vr > \Ov{\vr} - \frac{1}{N},     \end{array} \right.
\label{MHS}
\end{align}
where the constants $a_1$, $a_2$ are chosen in such a way that $p^N_{HS} \in C^1[0, \infty)$. The associated internal energy reads
\[
e_N(\vr, \vt) = e_m (\vr, \vt) + e_r(\vr, \vt) + \frac{1}{N} \vr + \int_{\Ov{\vr}/2}^\vr \frac{1}{z^2}\Big( ([ z - \Ov{\vr} - 1]^+)^\Gamma + p^N_{HS} (z) \Big) \D z
\]
cf. \cite[Chapter 8, Section 8.2]{FeiNovOpen}.

The weak formulation of the approximate problem reads:

	\begin{align}
	\intST{ \intO{ \left[ \vr \partial_t \varphi + \vr \vu \cdot \Grad \varphi \right] } } &= 0,
	\label{A3} \\
	\intST{ \intO{ \left[ b(\vr) \partial_t \varphi + b(\vr) \vu \cdot \Grad \varphi + \Big(
			b(\vr) - b'(\vr) \vr \Big) \Div \vu \varphi \right] } } &=0
	\label{A4}
\end{align}
for any $\varphi \in C^1(S_T \times \Ov{\Omega} )$, and any $b \in C^1(R)$, 	$b' \in C_c(R)$;

\begin{align}
	&\intST{ \intO{ \Big[ \vr \vu \cdot \partial_t \bfphi + \vr \vu \otimes \vu : \Grad \bfphi					 +
			p_N (\vr, \vt) \Div \bfphi \Big] } } \br &= \intST{ \intO{ \Big[ \mathbb{S}(\vt, \Ds \vu) : \Ds \bfphi - \vr \Grad G \cdot \bfphi \Big] } } - \frac{1}{\ep} \intST{ \int_{\partial \Omega}
		(\vu - \vuB) \cdot \bfphi \ \D \sigma_x }
	\label{A5}
\end{align}	
for any $\bfphi \in C^1(S_T \times \Ov{\Omega}; R^3)$, $\bfphi \cdot \vc{n}|_{\partial \Omega} = 0$;
\begin{align}
	- &\intST{ \intO{ \left[ \vr s(\vr, \vt) \partial_t \varphi + \vr s(\vr, \vt) \vu \cdot \Grad \varphi + \frac{\vc{q} (\vt, \Grad \vt)}{\vt} \cdot
			\Grad \varphi \right] } } \br &\geq \intST{ \intO{ \frac{\varphi}{\vt} \left[ \mathbb{S}(\vt, \Ds \vu) : \Ds \vu -
			\frac{\vc{q}(\vt, \Grad \vt) \cdot \Grad \vt }{\vt} \right] } }\br &+ \frac{1}{\ep}
	\intST{ \int_{\partial \Omega} \varphi \frac{|\vtB - \vt|^k (\vtB - \vt) }{\vt} \D \sigma_x }
	\label{A6}
\end{align}
for any $\varphi \in C^1(S_T \times \Ov{\Omega})$, $\varphi \geq 0$;
\begin{align}
	- &\intST{ \partial_t \psi	\intO{ \left[ \frac{1}{2} \vr |\vu|^2 + \vr e_N (\vr, \vt)  \right] } } + \frac{1}{\ep} \intST{ \psi  \int_{\partial \Omega} |\vt - \vtB|^k (\vt - \vtB)  \D \sigma_x
	}  \br
&+ \frac{1}{\ep} \intST{ \psi \int_{\partial \Omega} (\vu - \vuB) \cdot \vu \ \D \sigma_x }  \leq
	\intST{ \psi \intO{ \vr \vu \cdot \Grad G  } }
	\label{A7}
\end{align}
for any $\psi \in C^1(S_T)$, $\psi \geq 0$.

Relation \ref{A7} is the standard energy balance whereas the weak formulation is the same as in \cite{FeMuNoPo}
based on the abstract theory developed in \cite{FeNo6A}. The approximate system is ``almost closed'' in the sense that the energy flux through the boundary is controlled by the boundary conditions. Accordingly, solutions of the
approximate problem may be obtained exactly as in \cite{FeMuNoPo} as long as suitable {\it a priori} bounds are available. They are discussed in the forthcoming section.

\subsection{A priori bounds}

{\bf Step 1:}

For $\psi = 1$, the energy inequality \eqref{A7} yields
\begin{align}
\frac{1}{\ep} &\intST{ \int_{\partial \Omega} |\vt - \vtB|^k (\vt - \vtB)  \D \sigma_x } +
 \frac{1}{\ep} \intST{  \int_{\partial \Omega} (\vu - \vuB) \cdot \vu \ \D \sigma_x } \br &\leq
\intST{  \intO{ \vr \vu \cdot \Grad G  } } = 0	
\nonumber
\end{align}
as $G$ is independent of $t$. Consequently, since $\vt \geq 0$ we obtain
\begin{equation} \label{A8}
\| \vt \|_{L^{k + 1}(S_T \times \partial \Omega)} +
\| \vu \|_{L^2(S_T \times \partial \Omega; R^3)} \aleq 1,		
	\end{equation}
where the bound is independent of $\ep$, $N$. Here and hereafter, the symbol $A \aleq B$ means there is a positive
constant $c > 0$ such that $A \leq cB$.

\medskip

\noindent
{\bf Step 2:} In view of \eqref{A8}, the choice $\varphi = 1$ in the entropy inequality \eqref{A6} yields
\begin{equation} \label{A9}
\intST{ \intO{ \frac{1}{\vt} \left[ \mathbb{S}(\vt, \Ds \vu) : \Ds \vu -
		\frac{\vc{q}(\vt, \Grad \vt ) \cdot \Grad \vt }{\vt} \right] } } \leq c(\data, \ep),
\end{equation}		
where the right--hand side is independent of $N$ but may blow up for $\ep \to 0$.

\subsection{Existence of approximate solutions}

With \eqref{A8}, \eqref{A9} at hand, the remaining {\it a priori} bounds can be deduced exactly as in
\cite[Section 2.4]{FeMuNoPo}. Repeating step by step the arguments of the existence proof in \cite{FeMuNoPo}, we obtain the existence of approximate solutions.
\begin{Proposition} [{\bf Approximate solutions}] \label{AP1}
	
	In addition to the hypotheses of Theorem \ref{MT1}, let
	\[
		\ep > 0,\ \ N > 0
	\]
	be given.
	
	Then the approximate problem \eqref{A3} -- \eqref{A7} admits a solution $(\vr_{\ep,N}, \vt_{\ep, N}, \vu_{\ep, N})$.
	
	\end{Proposition}

\subsection{Limit $N \to \infty$}

With the {\it a priori} bounds \eqref{A8}, \eqref{A9} at hand, the limit $N \to \infty$ can be performed similarly to \cite[Chapter 8, Section 8.2]{FeiNovOpen}. We therefore obtain a family of approximate solutions
$(\vre, \vte, \vue)_{\ep > 0}$ satisfying the problem with the original pressure $p$ and the internal energy $e$.

\section{Limit $\ep \to 0$}
\label{e}

Our ultimate goal is to perform the limit $\ep \to 0$
in the family of approximate solutions $(\vre, \vte, \vue)$ obtained in the previous section.
They satisfy the following system of integral identities:

\begin{align}
	\intST{ \intO{ \left[ \vre \partial_t \varphi + \vre \vue \cdot \Grad \varphi \right] } } &= 0,
	\label{e1} \\
	\intST{ \intO{ \left[ b(\vre) \partial_t \varphi + b(\vre) \vue \cdot \Grad \varphi + \Big(
			b(\vre) - b'(\vre) \vre \Big) \Div \vue \varphi \right] } } &=0
	\label{e2}
\end{align}
for any $\varphi \in C^1(S_T \times \Ov{\Omega} )$, and any $b \in C^1(R)$, 	$b' \in C_c(R)$;
\begin{align}
	&\intST{ \intO{ \Big[ \vre \vue \cdot \partial_t \bfphi + \vre \vue \otimes \vue : \Grad \bfphi					 +
			p (\vre, \vte) \Div \bfphi \Big] } } \br &= \intST{ \intO{ \Big[ \mathbb{S}(\vte, \Ds \vue) : \Ds \bfphi - \vre \Grad G \cdot \bfphi \Big] } } - \frac{1}{\ep} \intST{ \int_{\partial \Omega}
		(\vue - \vuB) \cdot \bfphi \ \D \sigma_x }
	\label{e3}
\end{align}	
for any $\bfphi \in C^1(S_T \times \Ov{\Omega}; R^3)$, $\bfphi \cdot \vc{n}|_{\partial \Omega} = 0$;
\begin{align}
	- &\intST{ \intO{ \left[ \vre s(\vre, \vte) \partial_t \varphi + \vre s(\vre, \vte) \vue \cdot \Grad \varphi + \frac{\vc{q} (\vte, \Grad \vte)}{\vte} \cdot
			\Grad \varphi \right] } } \br &\geq \intST{ \intO{ \frac{\varphi}{\vte} \left[ \mathbb{S}(\vte, \Ds \vue) : \Ds \vue -
			\frac{\vc{q}(\vte, \Grad \vte) \cdot \Grad \vte }{\vte} \right] } }\br &+ \frac{1}{\ep}
	\intST{ \int_{\partial \Omega} \varphi \frac{|\vtB - \vte|^k (\vtB - \vte) }{\vte} \D \sigma_x }
	\label{e4}
\end{align}
for any $\varphi \in C^1(S_T \times \Ov{\Omega})$, $\varphi \geq 0$;
\begin{align}
	- &\intST{ \partial_t \psi	\intO{ \left[ \frac{1}{2} \vre |\vue|^2 + \vre e (\vre, \vte)  \right] } } + \frac{1}{\ep} \intST{ \psi  \int_{\partial \Omega} |\vte - \vtB|^k (\vte - \vtB)  \D \sigma_x
	}  \br
	&+ \frac{1}{\ep} \intST{ \psi \int_{\partial \Omega} (\vue - \vuB) \cdot \vue \ \D \sigma_x }  \leq
	\intST{ \psi \intO{ \vre \vue \cdot \Grad G  } }
	\label{e5}
\end{align}
for any $\psi \in C^1(S_T)$, $\psi \geq 0$.

In addition, the density is uniformly bounded,
\begin{equation} \label{BD}
	0 \leq \vre < \Ov{\vr} \ \mbox{a.a. in}\ (0,T) \times \Omega.
	\end{equation}

\subsection{Ballistic energy balance}

We start by choosing $\varphi = \psi(t) \tvt (t,x)$,
\[
\tvt \in C^1(S_T \times \Omega),\ \psi > 0,\  \tvt > 0,\ \tvt|_{\partial \Omega} = \vtB,
\]
as a test function in the approximate entropy balance \eqref{e4}. Adding the resulting expression to the energy
balance \eqref{e5} we obtain
\begin{align}
	- &\intST{ \partial_t \psi	\intO{ \left[ \frac{1}{2} \vre |\vue|^2 + \vre e - \tvt \vre s \right] } }
	+ \intST{ \psi
		\intO{ \frac{\tvt}{\vte}	 \left[ \mathbb{S}: \Ds \vue - \frac{\vc{q} \cdot \Grad \vte }{\vte} \right] } }  \br
	&+ \frac{1}{\ep} \intST{ \psi \int_{\partial \Omega} \frac{ |\vte - \vtB|^{k + 2} }{\vte} \ \D \sigma_x }
	+  \frac{1}{\ep} \intST{ \psi \int_{\partial \Omega} (\vue - \vuB) \cdot \vue \ \D \sigma_x }
	\br
	&\leq
	\intST{ \psi \intO{ \left[ \vre \vue \cdot \Grad G - \vre s \vue \cdot \Grad \tvt - \frac{\vc{q}}{\vte} \cdot \Grad \tvt
			- \partial_t \tvt \vre s
			\right] } }.
	\label{e6}
\end{align}

Next, the choice $\bfphi = \psi(t) \vuB(x)$ in the momentum balance \eqref{e3} yields
\begin{align}
	&\intST{ \intO{ \Big[ \vre \vue \cdot \vuB \partial_t \psi + \psi \vre \vue \otimes \vue : \Grad \vuB					 +
		\psi	p (\vre, \vte) \Div \vuB \Big] } } \br &= \intST{ \psi \intO{ \Big[ \mathbb{S}(\vte, \Ds \vue) : \Ds \vuB - \vre \Grad G \cdot \vuB \Big] } } - \frac{1}{\ep} \intST{ \psi \int_{\partial \Omega}
		(\vue - \vuB) \cdot \vuB \ \D \sigma_x }.
	\label{e7}
\end{align}
Note that $\bfphi = \psi(t) \vuB(x)$ is an admissible test function as $\vuB \cdot \vc{n}|_{\partial \Omega} = 0$.

Finally, taking $\varphi = \psi (t) \frac{1}{2} |\vuB|^2$ in the equation of continuity \eqref{e1} we get
\begin{equation} \label{e8}
- \intST{ \intO{ \left[ \frac{1}{2} \vre |\vuB|^2  \partial_t \psi + \psi \vre \vue \cdot \vuB \cdot \Grad \vuB \right] } } = 0.
\end{equation}

Summing up \eqref{e6}--\eqref{e8} we obtain an approximate ballistic energy balance in the form
\begin{align}
	- &\intST{ \partial_t \psi	\intO{ \left[ \frac{1}{2} \vre |\vue -\vuB|^2 + \vre e - \tvt \vre s \right] } }
	+ \intST{ \psi
		\intO{ \frac{\tvt}{\vte}	 \left[ \mathbb{S}: \Ds \vue - \frac{\vc{q} \cdot \Grad \vte }{\vte} \right] } }  \br
	&+ \frac{1}{\ep} \intST{ \psi \int_{\partial \Omega} \frac{ |\vte - \vtB|^{k + 2} }{\vte} \ \D \sigma_x }
	+  \frac{1}{\ep} \intST{ \psi \int_{\partial \Omega} |\vue - \vuB|^2  \ \D \sigma_x }
	\br
	&\leq
	\intST{ \psi \intO{ \left[ \vre (\vue - \vuB) \cdot \Grad G - \vre s \vue \cdot \Grad \tvt - \frac{\vc{q}}{\vte} \cdot \Grad \tvt
			- \partial_t \tvt \vre s
			\right] } } \br
	&- \intST{ \psi \intO{ \left[ \vre \vue \otimes \vue + p(\vre, \vte) \mathbb{I} - \mathbb{S}(\vte, \Ds
			\vue) \right] : \Ds \vuB } }\br &+ \frac{1}{2} \intST{ \psi \intO{\vre \vue \cdot \Grad |\vuB|^2 } }	
	\label{e9}
\end{align}
for any $\psi \in C^1(S_T)$, $\psi \geq 0$.	

\subsection{Extending boundary velocity}
\label{ext}

Following Galdi \cite[Lemma IX.4.1]{GALN}, see also
Kozono and Yanagisawa \cite[Proposition 1]{KozYan}, we may extend $\vuB$
\begin{equation} \label{e10}
	\vu_{\delta,B} (t,x) = {\bf curl}_x (d_\delta (x) \vc{z}_B(x)),\ {\bf curl_x}\  \vc{z}_B = \vuB,
\end{equation}
where
\begin{align}
	|d_\delta | &\leq 1,\ d_\delta (x) \equiv 1 \ \mbox{for all}\ x \ \mbox{in an open neighborhood of}\ \partial \Omega, \br
	d_\delta (x) &\equiv 0 \ \mbox{whenever}\ {\rm dist}[x, \partial \Omega] > \delta, \br
	|D^\alpha_x d_\delta (x) | &\leq c \frac{\delta}{{\rm dist}^{|\alpha|} [x, \partial \Omega]},\ |\alpha| = 1,2 ,\ 0 < \delta < 1,\
	x \in \Omega.
	\label{e11}	
\end{align}
The specific value of the parameter $\delta > 0$ will be fixed in the next section.

\subsection{Uniform bounds}

The desired uniform bounds will follow from the ballistic energy balance with the ansatz
\[
\psi = 1, \ \tvt = \vtB,
\]
where $\vtB$ is the harmonic extension of the boundary temperature,
\begin{equation} \label{e12}
\Del \vtB (t, \cdot) = 0 \ \mbox{for any}\ t \in S_T, \
\vtB \ \mbox{satisfies \eqref{dom} .}
\end{equation}

In accordance with the extension of $\vuB$ introduced in the preceding section, we have $\Div \vuB = 0$, and
the ballistic energy balance \eqref{e9} gives rise to
\begin{align}
	&\intST{
		\intO{ \frac{\vtB}{\vte}	 \left[ \mathbb{S} (\vte, \Ds \vue) : \Ds \vue + \frac{\kappa (\vte) |\Grad \vte|^2 }{\vte} \right] } }  \br
	&+ \frac{1}{\ep} \intST{ \int_{\partial \Omega} \frac{ |\vte - \vtB|^{\beta + 1} }{\vte} \ \D \sigma_x }
	+  \frac{1}{\ep} \intST{  \int_{\partial \Omega} |\vue - \vuB|^2  \ \D \sigma_x }
	\br
	&\leq -
	\intST{  \intO{ \left[ \vre \vuB \cdot \Grad G  + \vre s (\vre, \vte) \vue \cdot \Grad \vtB + \frac{\vc{q}}{\vte} \cdot \Grad \vtB
			+ \partial_t \vtB \vre s(\vre, \vte)
			\right] } } \br
	&- \intST{  \intO{ \left[ \vre \vue \otimes \vue  - \mathbb{S}(\vte, \Ds
			\vue) \right] : \Ds \vuB } },	
	\label{e13}
\end{align}	
where we have used
\[
\intST{ \intO{ \vre \vue \cdot \Grad G } } = 0 = \intST{ \intO{ \vre \vue \cdot \Grad |\vuB|^2 } }.
\]

Our goal is to show that all integrals on the right--hand side of \eqref{e12} are either bounded or dominated
be those on the left--hand side.

First, as a consequence of the standard maximum principle,
\[
\inf_{S_T \times \Omega} \vtB \geq \inf_{S_T \times \partial \Omega} \vtB > 0.
\]
Consequently, by virtue of hypothesis \eqref{m9} and the standard Sobolev embedding,
\begin{equation} \label{e14}
	\| \vue \|_{L^6(\Omega; R^3)}^2 \aleq \| \vue \|^2_{W^{1,2}(\Omega; R^3)} \aleq
\intO{ \frac{\vtB}{\vte} \mathbb{S} (\vte, \Ds \vue) : \Ds \vue    } + \int_{\partial \Omega} |\vue - \vuB|^2  \ \D \sigma_x + 1.
	\end{equation}
Similarly,
\begin{equation} \label{e15}
	\| \vte^{\frac{\beta}{2}} \|_{W^{1,2}(\Omega)}^2 \aleq \intO{ \frac{\vtB}{\vte^2} \kappa(\vte) |\Grad \vte |^2
	} +  \int_{\partial \Omega} \frac{ |\vte - \vtB|^{\beta + 1} }{\vte} \ \D \sigma_x + 1.
\end{equation}	

Next, it follows from the uniform bound \eqref{BD} that
\[
\left| \intST{  \intO{  \vre \vuB \cdot \Grad G }} \right| \aleq 1.
\]
In addition, by virtue of hypothesis \eqref{m9},
\[
\left| \intO{ \mathbb{S}(\vte, \Ds \vue) : \Ds \vuB } \right| \aleq \omega \| \Grad \vue \|^2_{L^2(\Omega; R^{3 \times 3})} + c(\omega) \| \vte \|^2_{L^2(\Omega)}
\]
for any $\omega > 0$. Choosing $\omega > 0$ small enough, this integral can be absorbed by the left--hand side of
\eqref{e13}.

Thus, in view of \eqref{e14}, \eqref{e15}, the previous observations, and boundedness of the density, inequality
\eqref{e13} gives rise to
\begin{align}
	&\intST{ \left( \|\vue \|^2_{W^{1,2}(\Omega)} + \| \vte^{\frac{\beta}{2}} \|_{W^{1,2}(\Omega)}^2 +
		\|\Grad \log (\vte) \|_{L^2(\Omega; R^3)}^2 \right)
	 }  \br
	&+ \frac{1}{\ep} \intST{ \int_{\partial \Omega} \frac{ |\vte - \vtB|^{\beta + 1} }{\vte} \ \D \sigma_x }
	+  \frac{1}{\ep} \intST{  \int_{\partial \Omega} |\vue - \vuB|^2  \ \D \sigma_x }
	\br
	&\aleq \left|
	\intST{  \intO{ \left[ \vre s (\vre, \vte) \vue \cdot \Grad \vtB + \frac{\vc{q}}{\vte} \cdot \Grad \vtB
			+ \partial_t \vtB \vre s(\vre, \vte)
			\right] } } \right| \br
	& + \left|  \intST{  \intO{ \vre (\vue \otimes \vue) : \Ds \vuB  } } \right| + 1.	
	\label{e17}
\end{align}	

\subsubsection{Entropy dependent terms}

Our goal is to control the first integral on the right--hand side of \eqref{e17}.
Denoting
\[
\mathcal{K}(\vt) = \int_1^\vt \frac{\kappa (z) }{z} \ \D z,
\]
we obtain,
\begin{align}
- \intO{ \frac{\vc{q}}{\vte} \cdot \Grad \vtB } &=
\intO{ \frac{\kappa (\vte) \Grad \vte }{\vte} \cdot \Grad \vtB } =
\intO{ \Grad \mathcal{K} (\vte) \cdot \Grad \vtB } \br &= \int_{\partial \Omega} \mathcal{K} (\vte) \Grad \vtB \cdot
\vc{n} \D \sigma_x,
\nonumber
\end{align}
where we have used the fact that $\vtB$ is harmonic inside $\Omega$. It follows from hypothesis \eqref{m9}
\[
\left| \int_{\partial \Omega} \mathcal{K} (\vte) \Grad \vtB \cdot
\vc{n} \D \sigma_x \right| \aleq \left(1+ \intST{  \int_{\partial \Omega} \frac{ |\vte - \vtB|^{\beta + 1} }{\vte} \ \D \sigma_x}\right).
\]
Consequently, inequality \eqref{e17} reduces to
\begin{align}
	&\intST{ \left( \|\vue \|^2_{W^{1,2}(\Omega)} + \| \vte^{\frac{\beta}{2}} \|_{W^{1,2}(\Omega)}^2 +
		\|\Grad \log (\vte) \|_{L^2(\Omega; R^3)}^2 \right)
	}  \br
	&+ \frac{1}{\ep} \intST{ \int_{\partial \Omega} \frac{ |\vte - \vtB|^{\beta + 1} }{\vte} \ \D \sigma_x }
	+  \frac{1}{\ep} \intST{  \int_{\partial \Omega} |\vue - \vuB|^2  \ \D \sigma_x }
	\br
	&\aleq \left|
	\intST{  \intO{ \left[ \vre s (\vre, \vte) \vue \cdot \Grad \vtB
			+ \partial_t \vtB \vre s(\vre, \vte)
			\right] } } \right| \br
	& + \left|  \intST{  \intO{ \vre (\vue \otimes \vue) : \Ds \vuB  } } \right| + 1.	
	\label{e18}
\end{align}

Next, using hypothesis \eqref{m8}, we easily obtain
\begin{equation} \label{e19}
	0 \leq \vr s_m (\vr, \vt) \aleq \left( \vr \log^+(\vr) + \vr \log^+ (\vt) + 1 \right),
	\end{equation}
see \cite[Chapter 4, Section 12.4.2]{FeiNovOpen}. Consequently, by means of \eqref{BD},
\[
|\vre s(\vre, \vte) \vue| \aleq |\vre s_m (\vre, \vte) \vue | + |\vte|^3 |\vue| \aleq
(1 + |\vte|^3) |\vue|.
\]
Since $\beta > 6$, the integral
\[
\intST{ \intO{ (1 + |\vte|^3) |\vue|}}
\]
can be absorbed by the left--hand side of \eqref{e18}. The integral
\[
\intST{ \intO{ \partial_t \vtB \vre s(\vre, \vte) }}
\]
can be handled in a similar fashion. We conclude
\begin{align}
	&\intST{ \left( \|\vue \|^2_{W^{1,2}(\Omega)} + \| \vte^{\frac{\beta}{2}} \|_{W^{1,2}(\Omega)}^2 +
		\|\Grad \log (\vte) \|_{L^2(\Omega; R^3)}^2 \right)
	}  \br
	&+ \frac{1}{\ep} \intST{ \int_{\partial \Omega} \frac{ |\vte - \vtB|^{\beta + 1} }{\vte} \ \D \sigma_x }
	+  \frac{1}{\ep} \intST{  \int_{\partial \Omega} |\vue - \vuB|^2  \ \D \sigma_x }
	\br& \aleq \left| \intST{  \intO{ \vre (\vue \otimes \vue) : \vuB  } } \right| + 1.	
	\label{e19a}
\end{align}

\subsubsection{Convective term}
\label{ccc}

It remains to handle the convective term
\[
\vre \vue \otimes \vue = \vre (\vue - \vuB) \otimes (\vue - \vuB) + \vre \vuB \otimes \vue
+ \vre (\vue - \vuB) \otimes \vuB.
\]
Writing
\begin{align}
\intO{ \vre (\vue \otimes \vue) : \Ds \vuB  }  &=
\intO{ \vre (\vue - \vuB) \otimes (\vue - \vuB) : \Ds \vuB  } \br
& +\intO{ \vre \vuB \otimes \vue : \Ds \vuB  } + \intO{ \vre (\vue - \vuB) \otimes \vuB : \Ds \vuB  }
\nonumber
\end{align}
we immediately see that the only problematic term is
\[
\intO{ \vre (\vue - \vuB) \otimes (\vue - \vuB) : \Ds \vuB  }.
\]

Let $\vc{w}_\ep$ be the unique solution of the Dirichlet problem
\begin{equation} \label{e100}
	\Del \vc{w}_\ep = 0,\ \vc{w}_{\ep}|_{\partial \Omega} = (\vue - \vuB).
	\end{equation}
Write
\begin{align}
&\intO{ \vre (\vue - \vuB) \otimes (\vue - \vuB) : \Ds \vuB  } =
\intO{ \vre (\vue - \vuB - \vc{w}_\ep) \otimes (\vue - \vuB - \vc{w}_\ep ) : \Ds \vuB  } \br
&\quad + \intO{ \vre \vc{w}_\ep \otimes (\vue - \vuB) : \Ds \vuB }
+ \intO{ \vre (\vue - \vuB - \vc{w}_\ep ) \otimes \vc{w}_\ep : \Ds \vuB }.
\label{e101}
\end{align}
By means of the uniform bound
on the density, we get
\begin{align}
&\left| \intO{ \vre (\vue - \vuB - \vc{w}_\ep) \otimes (\vue - \vuB - \vc{w}_\ep ) : \Ds \vuB  }  \right| \br
&\quad \aleq \intO{ \frac{|\vue - \vuB - \vc{w}_\ep|^2}{ {\rm dist}^2(x, \partial \Omega) } {\rm dist}^2(x, \partial \Omega) | \Ds \vuB |  }.
\nonumber
\end{align}
By virtue of Hardy--Sobolev inequality, we obtain
\[
 \intO{ \frac{|\vue - \vuB - \vc{w}_\ep|^2}{ {\rm dist}^2(x, \partial \Omega) } } \aleq
 \| \vue - \vuB - \vc{w}_\ep \|_{W^{1,2}_0(\Omega; R^3)}^2 	\aleq \| \vue \|^2_{W^{1,2}(\Omega; R^3)} + 1.
\]
Thus going back to Section \ref{ext}, we can choose $\delta = \delta(\omega) > 0$ so small that
\begin{equation} \label{e102}
\left| \intO{ \vre (\vue - \vuB - \vc{w}_\ep) \otimes (\vue - \vuB - \vc{w}_\ep ) : \Ds \vuB  }  \right|
\leq \delta \| \vue \|^2_{W^{1,2}(\Omega; R^3)} + c.
	\end{equation}
for any $\omega > 0$.

To control the remaining two integrals in \eqref{e101}, we first evoke the $L^p-$elliptic estimates
applied to \eqref{e100}:
\[
\| \vc{w}_\ep \|_{W^{1,p}(\Omega)} \aleq \| \vue - \vuB \|_{W^{1 - \frac{1}{p}, p}(\partial \Omega)}.
\]
In particular, for $p = \frac{4}{3}$, we get
\[
\| \vc{w}_\ep \|_{W^{1,\frac{4}{3}}(\Omega)} \aleq \| \vue - \vuB \|_{W^{\frac{1}{4}, \frac{4}{3}}(\partial \Omega)}.
\]
Moreover, we use the interpolation to obtain the following inequality 
\[
\|  \vue - \vuB  \|_{W^{\frac{1}{4}, \frac{4}{3} }(\partial \Omega)} \aleq
\|  \vue - \vuB  \|_{W^{\frac{1}{4}, 2 }(\partial \Omega)} \leq \|  \vue - \vuB  \|_{W^{\frac{1}{2}, 2 }(\partial \Omega)}^{\frac{1}{2}} \|  \vue - \vuB  \|_{L^2(\partial \Omega)}^{\frac{1}{2}}.
\]
Finally, we recall the standard Sobolev embedding
\[
\| \vc{w}_\ep \|_{L^q(\Omega)} \aleq \| \vc{w}_\ep \|_{W^{1,\frac{4}{3}}(\Omega)},\
1 \leq q \leq  = \frac{12}{5}.
\]
Thus we may infer that
\begin{align}
	\| \vc{w}_\ep \|_{L^2(\Omega, R^3)} &\aleq  \|  \vue - \vuB  \|_{W^{\frac{1}{2}, 2 }(\partial \Omega; R^3)}^{\frac{1}{2}} \|  \vue - \vuB  \|_{L^2(\partial \Omega; R^3)}^{\frac{1}{2}}  \br &\aleq
\|  \vue - \vuB  \|_{W^{1, 2 }(\Omega; R^3)}^{\frac{1}{2}} \|  \vue - \vuB  \|_{L^2(\partial \Omega; R^3)}^{\frac{1}{2}}
\label{e103}
\end{align}
where we have here used the trace theorem.

Going back to \eqref{e101}, we have to estimate the products
\[
\intO{ |\vue| |\vc{w}_\ep |} \leq \| \vue \|_{L^2(\Omega; R^3)}  \| \vc{w}_\ep \|_{L^2(\Omega; R^3)}.
\]
It follows from \eqref{e103}
\begin{align}
\| \vue \|_{L^2(\Omega; R^3)} & \| \vc{w}_\ep \|_{L^2(\Omega; R^3)} \aleq
\| \vue \|_{L^2(\Omega; R^3)} \|  \vue - \vuB  \|_{W^{1, 2 }(\Omega; R^3)}^{\frac{1}{2}} \|  \vue - \vuB  \|_{L^2(\partial \Omega; R^3)}^{\frac{1}{2}} \br
&\leq \omega \| \vue \|_{L^2(\Omega; R^3)} \|  \vue - \vuB  \|_{W^{1, 2 }(\Omega; R^3)} +
c(\omega) \| \vue \|_{L^2(\Omega; R^3)}  \|  \vue - \vuB  \|_{L^2(\partial \Omega; R^3)},
\nonumber	
	\end{align}
where the first term on the right--hand side may be absorbed by the left--hand side of \eqref{e19a} if $\omega > 0$ is small enough. Finally, repeating the same argument,
\[
 \| \vue \|_{L^2(\Omega; R^3)}  \|  \vue - \vuB  \|_{L^2(\partial \Omega; R^3)} \leq
 \omega  \| \vue \|_{L^2(\Omega; R^3)}^2 + c(\omega)\|  \vue - \vuB  \|_{L^2(\partial \Omega; R^3)}^2
 \]
we conclude that also this term is controlled by the left--hand side of \eqref{e19a} if $\omega > 0$ is fixed
small enough and $\ep \to 0$.

In view of the preceding argument, inequality \eqref{e19a} gives rise to the desired conclusion
\begin{align}
	&\intST{ \left( \|\vue \|^2_{W^{1,2}(\Omega)} + \| \vte^{\frac{\beta}{2}} \|_{W^{1,2}(\Omega)}^2 +
		\|\Grad \log (\vte) \|_{L^2(\Omega; R^3)}^2 \right)
	}  \br
	&+ \frac{1}{\ep} \intST{ \int_{\partial \Omega} \frac{ |\vte - \vtB|^{\beta + 1} }{\vte} \ \D \sigma_x }
	+  \frac{1}{\ep} \intST{  \int_{\partial \Omega} |\vue - \vuB|^2  \ \D \sigma_x } \aleq 1	
	\label{e20}
\end{align}
uniformly for $\ep \to 0$.

\subsection{Pressure estimates}

The uniform bounds \eqref{BD}, \eqref{e20} are strong enough to control all terms in the field equations with the only exception of the pressure and the associated integral energy that are singular for $\vr \to \Ov{\vr}$.
The desired estimates can be derived following the arguments of \cite{BreFeiNov20}.

First, we introduce the so--called Bogovskii operator
\begin{align}
	\mathcal{B} : L^q_0 (\Omega) &\equiv \left\{ f \in L^q(\Omega) \ \Big| \ \intO{ f } = 0 \right\}
	\to W^{1,q}_0(\Omega, R^d),\ 1 < q < \infty,\br
	\Div \mathcal{B}[f] &= f .
	\nonumber
\end{align}
The operator $\mathcal{B}$ can be constructed by means of the original ansatz of Bogovskii \cite{BOG}
elaborated by Galdi \cite[Chapter 3]{GALN}, and later revisited by Geissert, Heck, and Hieber \cite{GEHEHI}.
$\mathcal{B}$ maps bounded sets of $L^q_0(\Omega)$ into bounded sets of $W^{1,q}_0(\Omega; R^d)$ and
bounded sets in $(W^{1,q}(\Omega))^* \perp 1$ to bounded sets of $L^q(\Omega; R^d)$ provided $\Omega$ is a Lipschitz domain.

\subsubsection{Integrability of the pressure}

We consider the quantity
\[
\bfphi = \mathcal{B} \left[ \vre - \frac{1}{|\Omega|} \intO{\vre } \right]
\]
as a test function in the approximate momentum balance \eqref{e3}. After a straightforward manipulation, we get
\begin{align}
	&\intST{ \intO{ p(\vre, \vte) \left[ \vre - \frac{1}{|\Omega|} \intO{\vre } \right] } }  \br
	&= \intST{ \intO{ \Big[ \vre \vue \cdot \mathcal{B}[ \Div (\vre \vue) ] - \vre \vue \otimes \vue : \Grad \mathcal{B} \left[ \vre - \frac{1}{|\Omega|} \intO{\vre } \right] 				
			 \Big] } } \br &\quad + \intST{ \intO{  \mathbb{S}(\vte, \Ds \vue) : \Ds \mathcal{B} \left[ \vre - \frac{1}{|\Omega|} \intO{\vre } \right] } } \br &\quad  - \intST{ \intO{ \vre \Grad G \cdot \mathcal{B} \left[ \vre - \frac{1}{|\Omega|} \intO{\vre } \right]   } }
	\label{e22}
\end{align}	
As $\vre$ are uniformly bounded, it is easy to check that the right--hand side of \eqref{e22} is bounded
by means of the uniform estimates established in \eqref{e20}. Since the total mass is constant, we get
\[
\frac{1}{|\Omega|} \intO{ \vre } < \Ov{\vr}.
\]
In particular, boundedness of the integral on the right--hand side of \eqref{e22} yields a uniform bound
\begin{equation} \label{e23}
\intST{ \intO{ p(\vre, \vte) }} \aleq 1	
	\end{equation}
uniformly for $\ep \to 0$. Moreover, as observed in \cite[Section 2.1, formula (2.4)]{BreFeiNov20},
uniform integrability of the hard--sphere pressure yields equi--integrability of the associated internal energy.
Specifically,
\begin{equation} \label{e24}
(\vre e(\vre, \vte) )_{\ep > 0} \ \mbox{is}\ L^1 - \mbox{equi--integrable.}	
	\end{equation}

Finally, using \eqref{e24}, we may deduce from the ballistic energy inequality \eqref{e9}, exactly as for the
initial--value problem,
\begin{equation} \label{e27}
	(\vre |\vue|^2 )_{\ep > 0} , \
	(\vre e(\vre, \vte))_{\ep > 0} \ \mbox{bounded in}\ L^\infty(S_T; L^1(\Omega)).
\end{equation}

\subsubsection{Equi--integrability of the pressure}

The ultimate goal of this section is to establish $L^1$--integrability of the pressure. To this end, we repeat the
above procedure with the test function
\[
\bfphi = \mathcal{B} \left[ b(\vre) - \frac{1}{|\Omega|} \intO{ b(\vre) } \right]
\]
in \eqref{e3}, where $b$ is a function compatible with the renormalized equation of continuity \eqref{e2}.
After a straightforward manipulation, we get
\begin{align}
	&\intST{ \intO{ p(\vre, \vte) \left[  b(\vre) - \frac{1}{|\Omega|} \intO{ b(\vre) } \right] } }  \br
	&= - \intST{ \intO{ \vre \vue \otimes \vue : \Grad \mathcal{B} \left[ b(\vre) - \frac{1}{|\Omega|} \intO{ b(\vre) } \right] 				
			 } } \br &\quad + \intST{ \intO{  \mathbb{S}(\vte, \Ds \vue) : \Ds \mathcal{B} \left[ b(\vre) - \frac{1}{|\Omega|} \intO{ b(\vre) } \right] } } \br &\quad  - \intST{ \intO{ b(\vre) \Grad G \cdot \mathcal{B} \left[ b(\vre) - \frac{1}{|\Omega|} \intO{ b(\vre) } \right]   } }\br
		 &\quad+ \intST{ \intO{ \vre \vue  \cdot \mathcal{B} \left[ ( b(\vre) - b'(\vre) \vre) \Div \vue
		 		- \frac{1}{|\Omega|} \intO{(b(\vre) - b'(\vre) \vre) \Div \vue } \right] } } \br
	 	&\quad - \intST{ \intO{ \vre \vue \cdot \mathcal{B}[ \Div (b(\vre) \vue) ]     }}.
	\label{e25}
\end{align}	

The main idea is to consider $b(\vr) = p_{HS}^\nu (\vr)$, where $\nu > 0$ is sufficiently small. In view of \eqref{e23}, $p_{HS}^\nu (\vre)$ is uniformly bounded in the Lebesgue space $L^q((0,T) \times \Omega)$, where $q > 1$ can be arbitrarily large provided $\nu > 0$ is small. It would follow that
\[
p(\vre, \vte)^{1 + \nu} \ \mbox{is bounded in}\ L^1((0,T) \times \Omega)
\]
provided we can show that all integrals on the right--hand side of \eqref{e25} are bounded uniformly for $\ep \to 0$.

There are certain technical difficulties to carry out the above programme. To begin, $p_{HS}$ is singular therefore
not directly eligible for the renormalized equation \eqref{e2}. Fortunately, we may consider $b$ any $C^1$ truncation of $p_{HS}$ as $\vre$ are uniformly bounded and behaviour of $b$ for large argument is irrelevant.
In particular, we may consider
\[
b(\vre) = ( p^N_{HS} )^\nu(\vre),
\]
with $P^N_{HS}$ given by \eqref{MHS}. The desired conclusion then follows as soon as we are able to show that
the right--hand side of \eqref{e25} remains bounded uniformly for $N \to \infty$, $\ep \to 0$. Given the
properties of the operator $\mathcal{B}$, the only problematic term is the integral
\begin{align}
\int_{S_T} &\int_{\Omega} \vre \vue  \cdot \mathcal{B} \Big[ ( (p_{HS}^N)^\nu (\vre) - [(p_{HS}^N)^\nu]'(\vre) \vre) \Div \vue \br
		&- \frac{1}{|\Omega|} \intO{((p_{HS}^N)^\nu(\vre) - [ (p_{HS}^N)^\nu ]'(\vre) \vre) \Div \vue } \Big] \dxdt.
\label{int}
\end{align}
Indeed this integral contains the derivative $[ (p_{HS}^N)^\nu ]'$ that is more singular than
$(p_{HS}^N)^\nu$ in the neighbourhood of $\Ov{\vr}$ when $N \to \infty$.

On the one hand,
\begin{equation} \label{NHP}
 [ (p_{HS}^N)^\nu ]'(\vr) \approx \left\{ \begin{array}{l} (\Ov{\vr} - \vr)^{-(\nu s + 1)} \ \mbox{if} \
 	\vr \leq \Ov{\vr} - \frac{1}{N} \\ \\ N^{\nu s + 1}, \ \mbox{if} \
 	\vr > \Ov{\vr} - \frac{1}{N}, \end{array} \right.
\end{equation}
On the other hand, in accordance with \eqref{e27},
\[
	(\vre e(\vre, \vte))_{\ep > 0} \ \mbox{bounded in}\ L^\infty(S_T; L^1(\Omega)).
	\]
Thus, by means of hypothesis \eqref{m3},
\begin{equation} \label{e28}
	\left( (\Ov{\vr} - \vr)^{1 - s} \right)_{\ep > 0} \ \mbox{bounded in}\
	L^\infty(S_T; L^1(\Omega)).
	\end{equation}
Comparing \eqref{NHP} with \eqref{e28} we conclude
\begin{equation} \label{e33}
[(p_{HS}^N)^\nu]'(\vre) \vre) \ \mbox{bounded in}\ L^\infty(S_T; L^q(\Omega))
\ \mbox{provided}\ s \geq q \nu s + q + 1.
\end{equation}	

If $s > 3$, we can find $\nu$ small enough so that $q > 2$ in \eqref{e33}. Consequently,
\[
( [(p_{HS}^N)^\nu]'(\vre) \vre) \Div \vue )_{\ep > 0} \ \mbox{is bounded in}\
L^2(S_T; L^\beta(\Omega)) \ \mbox{for some} \ \beta > 1;
\]
whence
\[
\mathcal{B} \Big[ ( (p_{HS}^N)^\nu (\vre) - [(p_{HS}^N)^\nu]'(\vre) \vre) \Div \vue
- \frac{1}{|\Omega|} \intO{((p_{HS}^N)^\nu(\vre) - [ (p_{HS}^N)^\nu ]'(\vre) \vre) \Div \vue } \Big]
\]
is bounded in $L^2(S_T; L^{\frac{3}{2}}(\Omega; R^3))$. In view of the uniform bounds \eqref{e20},
we conclude that the integral \eqref{int} remains bounded uniformly for $N \to \infty$, $\ep \to 0$.
Thus we have obtained the desired conclusion
\begin{equation} \label{e36}
\intST{ \intO{ p(\vre, \vte)^{1 + \nu} } } \aleq 1 \ \mbox{for some}\ \nu > 0. 	
	\end{equation}

\subsection{Conclusion, limit $\ep \to 0$}

In the previous section, we have obtained all available uniform bounds, specifically, the uniform density estimates \eqref{BD}, the ``dissipative estimates'' \eqref{e20}, and the pressure estimates \eqref{e36}.
Now it is standard, extracting suitable subsequences as the case may be, to identify the limits
\begin{align}
	\vre &\to \vr \ \mbox{in}\ C_{\rm weak}(S_T; L^q(\Omega)) \ \mbox{for any finite}\ q, \br
	\vue &\to \vu \ \mbox{weakly in}\ L^2(0,T; W^{1,2}(\Omega; R^3) , \br
	\vte &\to \vt \ \mbox{weakly in}\ L^2(0,T; W^{1,2}(\Omega) .
	\label{e37}
	\end{align}
Moreover, it follows from \eqref{e20} that the limit velocity and temperature satisfy the desired
boundary conditions
\begin{equation} \label{e38}
	\vu|_{\partial \Omega} = \vuB,\ \vt|_{\partial \Omega} = \vtB.
	\end{equation}

Finally, it is a routine matter to perform the limit in the approximate equation of continuity
\eqref{e1}, \eqref{e2}, the momentum equation \eqref{e3}, the entropy inequality \eqref{e4}, and the
ballistic energy balance \eqref{e9} as long as we can show strong (pointwise a.a.) convergence
\begin{equation} \label{e39}
	\vre \to \vr, \ \vte \to \vt \ \mbox{a.a. in}\ (0,T) \times \Omega.
	\end{equation}
This is a non--trivial task, however nowadays well understood. In particular, the compactness arguments based on Div-curl Lemma and Lions' identity can be modified to accommodate the time periodic
setting exactly as in \cite[Section 9.3]{FeMuNoPo}.

We have proved Theorem \ref{MT1}.

\section{Concluding remarks}
\label{CR}

As already pointed out, the above result can be easily extended to the case of time periodic boundary velocity as well as time periodic potential $G$. A more delicate issue would be considering general inflow/outflow boundary
conditions in the spirit of \cite{AbbFei21per}. The present approximation does not apply as the total mass is generally not conserved. Adopting the indirect method of \cite{AbbFei21per} based applying a fixed point argument to the associated Poincar\' e map would result in essential technical difficulties due to the presence of the
internal energy equation.

The possibility of eliminating the hard--sphere pressure component remains largely open as long as the boundary
velocity is non--zero. The main and possibly the only unsurmountable problem is controlling the convective term
as in Section \ref{ccc}. Note that the same difficulty arises at the level of stationary solutions, cf.
\cite{CiFeJaPe1}.



\def\cprime{$'$} \def\ocirc#1{\ifmmode\setbox0=\hbox{$#1$}\dimen0=\ht0
	\advance\dimen0 by1pt\rlap{\hbox to\wd0{\hss\raise\dimen0
			\hbox{\hskip.2em$\scriptscriptstyle\circ$}\hss}}#1\else {\accent"17 #1}\fi}

\end{document}